\documentclass[a4paper]{amsart}

\newtheorem{e-proposition}[theorem]{Proposition}

\newtheorem{e-definition}[theorem]{Definition\rm}

\newtheorem{theoreme}{Th\'eor\`eme}[section]
\newtheorem{lemme}[theoreme]{Lemme}
\newtheorem{proposition}[theoreme]{Proposition}

\newcommand{\taum}{\tau^{-}}
\newcommand{\clWalc}{\overline{\mathcal{W}}}
\newcommand{\Walc}{\mathcal{W}}
\newcommand{\w}{\omega}

\begin{document}

\title[Un th\'eor\`eme de convexit\'e r\'eel pour un moment quasi-hamiltonien]{Un th\'eor\`eme de convexit\'e r\'eel pour\\ les applications moment \`a valeurs\\ dans un groupe de Lie}

\author{Florent Schaffhauser}

\address{Dept. of Mathematics - Keio University, Hiyoshi 3-14-1,
Kohoku-ku 223-8522, Yokohama, Japan}
\email{florent@math.jussieu.fr}

\thanks{Supported by the Japanese Society for Promotion of Science (JSPS)}

\begin{abstract}
Dans cette note, on pr\'esente l'\'enonc\'e et les principales id\'ees de la d\'emonstration d'un th\'eor\`eme de convexit\'e r\'eel pour les applications moment \`a valeurs dans un groupe de Lie. Ce r\'esultat est un analogue quasi-hamiltonien du th\'eor\`eme de O'Shea et Sjamaar dans le cadre hamiltonien usuel. On d\'emontre que l'image par l'application moment du lieu des points fixes d'une involution renversant la $2$-forme de structure d'un espace quasi-hamiltonien est un polytope convexe, et l'on d\'ecrit ce polytope comme sous-polytope du polytope moment.
\vskip 0.5\baselineskip

In this note, we state and give the main ideas of the proof of a real convexity theorem for group-valued momentum maps. This result is a quasi-Hamiltonian analogue of the O'Shea-Sjamaar theorem in the usual Hamiltonian setting. We prove here that the image under the momentum map of the fixed-point set of a form-reversing involution defined on a quasi-Hamiltonian space is a convex polytope, that we describe as a subpolytope of the momentum polytope.
\end{abstract}

\maketitle

\section*{Abridged English version}

Let $(G,(.\ |\ .))$ be a compact connected Lie group whose Lie algebra $\mathfrak{g}:=Lie(G)$ is endowed with an $Ad$-invariant (positive definite) scalar product $(.\ |\ .)$. We denote by $\theta^L=g^{-1}.dg$ and by $\theta^R=dg.g^{-1}$ the Maurer-Cartan $1$-forms on $G$ and by $\chi$ the Cartan $3$-form, that is to say the bi-invariant $3$-form defined for $X,Y,Z\in\mathfrak{g}=T_1 G$ by $\chi_1(X,Y,Z):=([X,Y]\ |\ Z)$. When $G$ acts on a manifold $M$, we denote by $X^{\sharp}_x:=\frac{d}{dt}|_{t=0}(\exp(tX).x)$ the fundamental vector field associated to $X\in\mathfrak{g}$ by the action of $G$ on $M$. A quasi-Hamiltonian $G$-space $(M,\w,\mu:M\to U)$ is a manifold $M$, acted upon by the group $G$, endowed with an invariant $2$-form $\w$ such that there exists an equivariant map $\mu:M\to G$ (for the conjugacy action of $G$ on itself) satisfying~:
\begin{enumerate}
\item[(i)] $d\w=-\mu^*\chi$
\item[(ii)] for all $x\in M$, $\ker \w_x =\{X^{\sharp}_x\ :\ X\in\mathfrak{g}\, |\, Ad\, \mu(x).X=-X\}$
\item[(iii)] for all $X \in\mathfrak{g}$, the interior product $\iota_{X^{\sharp}}\w$ satisfies $\iota_{X^{\sharp}}\w = \frac{1}{2}\mu^*(\theta^L + \theta^R\ |\ X)$
\end{enumerate}
This notion was introduced by Alekseev, Malkin and Meinrenken in \cite{AMM} and basic examples include the conjugacy classes of the Lie group $G$. The map $\mu$ is called the \emph{momentum map}, underlining the important analogies with usual Hamiltonian spaces that one encounters in the study of quasi-Hamiltonian spaces. In the rest of this paper, the compact connected Lie group $G$ will be assumed to be simply connected as well. In \cite{AMM}, Alekseev, Malkin and Meinrenken showed that in this case the intersection of the image $\mu(M)$ of the momentum map $\mu$ with a fundamental domain $\exp(\clWalc)$ for the conjugacy action ($\Walc\subset \mathfrak{g}$ is a Weyl alcove in the Lie algebra of $G$) is a convex polytope (a result which they derived from Meinrenken and Woodward in \cite{MW2}). Here, we study the image $\mu(M^{\beta})$ under the momentum map of the fixed-point set $M^{\beta}$ of an involution $\beta$ defined on $M$ satisfying $\beta^*\w=-\w$ (we shall say that $\beta$ reverses the structural $2$-form $\w$) and additional compatibility relations with the action and the momentum map (see theorem \ref{convexite_max}). To formulate these compatibility conditions, the group $G$ is assumed to be endowed with an involutive automorphism $\tau$. Additionally, we assume that the involution $\taum(g):=\tau(g^{-1})$ leaves a \emph{maximal} torus $T$ of $G$ \emph{pointwise} fixed. Such an involution always exists on a compact connected simply connected Lie group (see \cite{Lo}). Under this assumption, we show (theorem \ref{convexite_max}) that the set $\mu(M^\beta)\cap\exp(\clWalc)$ is a convex set, and that it is in fact equal to the whole momentum polytope $\mu(M)\cap\exp(\clWalc)$. The proof is sketched in section \ref{preuve} and details are available in \cite{these_Sch}. It rests on the existence of a symplectic cross-section $N\subset M$, to which one applies Duistermaat theorem in \cite{Du} to show that~: $$\mu(M^{\beta})\cap\exp(\clWalc)=\overline{\mu(N^{\beta})}=\overline{\mu(N)}=\mu(M)\cap\exp(\clWalc)$$ In a forthcoming paper, we shall try to show that when we drop the assumption that $\taum$ leaves pointwise fixed a \emph{maximal} torus $T$ of $G$, we still have that $\mu(M^{\beta})\cap\exp(\clWalc)$ is a convex subpolytope of the momentum polytope $\mu(M)\cap\exp(\clWalc)$, obtained, in an exact analogy with the usual Hamiltonian case treated in \cite{OSS}, by intersecting it with $Fix(\taum)$. Applications of theorem \ref{convexite_max} include the construction of Lagrangian submanifolds in representation spaces of surface groups, for which we refer to \cite{these_Sch}, where we give an example of an involution $\beta$ satisfying the assumptions of the theorem.

\section{Introduction}

Les propri\'et\'es de convexit\'e de l'image d'une application moment, qu'elle soit \`a valeurs dans le dual d'une alg\`ebre de Lie (cas hamiltonien usuel) ou \`a valeurs dans un groupe de Lie (cas quasi-hamiltonien) ont \'et\'e \'etudi\'ees de divers points de vue depuis le d\'ebut des ann\'ees 1980. Le premier r\'esultat a \'et\'e obtenu ind\'ependamment par Atiyah dans \cite{Ati} et par Guillemin et Sternberg dans \cite{GS1}, et concerne les actions hamiltoniennes de groupes compacts connexes \emph{ab\'eliens}, c'est-\`a-dire des actions hamiltoniennes de tores, sur des vari\'et\'es symplectiques connexes compactes. Le cas des actions hamiltoniennes de tores sur des vari\'et\'es symplectiques connexes non n\'ecessairement compactes (l'application moment est alors suppos\'ee propre) a \'et\'e trait\'e plus tard par Condevaux, Dazord et Molino dans \cite{CDM}; leurs id\'ees ont \'et\'e reprises et d\'evelopp\'ees par Hilgert, Neeb et Plank dans \cite{HNP}. Le cas des actions hamiltoniennes de groupes compacts connexes \emph{non-ab\'eliens} a \'et\'e trait\'e par Kirwan dans \cite{Kir}, prouvant une conjecture de Guillemin et Sternberg dans \cite{GS2} (o\`u est trait\'e le cas k\"ahl\'erien). Le th\'eor\`eme de Kirwan a ensuite \'et\'e \'etendu aux vari\'et\'es symplectiques connexes non n\'ecessairement compactes dans \cite{CDM,HNP} et \cite{Sjamaar}. Si l'on s'int\'eresse maintenant aux actions \emph{quasi}-hamiltoniennes de groupes compacts connexes, il faut, pour pouvoir parler de convexit\'e, supposer de plus que le groupe consid\'er\'e est \emph{simplement connexe}. En effet, le cas hamiltonien montre que les propri\'et\'es de convexit\'e de l'image de l'application moment sont obtenues en consid\'erant l'intersection de cette image avec un domaine fondamental pour l'action co-adjointe (c'est-\`a-dire l'action du groupe sur l'espace o\`u l'application moment prend ses valeurs), y compris dans le cas ab\'elien, o\`u l'action co-adjointe est triviale. En proc\'edant de m\^eme dans le cas quasi-hamiltonien, c'est-\`a-dire en consid\'erant l'intersection de l'image du moment avec un domaine fondamental de l'action du groupe sur lui-m\^eme par conjugaison, on constate qu'un tel domaine fondamental s'identifie \`a une partie convexe de l'alg\`ebre de Lie du groupe consid\'er\'e si et seulement si celui-ci est simplement connexe (voir par exemple \cite{Bourbaki} ou \cite{Lo}). En cons\'equence, pour que la notion de convexit\'e ait un sens pour les applications moment \`a valeurs dans un groupe compact, on se limitera aux groupes compacts connexes et simplement connexes.\\
\indent Fixons un groupe $G$ compact connexe et simplement connexe. Pour un tel groupe, le cas hamiltonien usuel distingue deux types de r\'esultats~: le th\'eor\`eme de convexit\'e de Guillemin-Sternberg-Kirwan et le th\'eor\`eme de convexit\'e dit \emph{r\'eel} de O'Shea et Sjamaar (\cite{OSS}) o\`u l'on \'etudie l'image par l'application moment du lieu des points fixes d'une involution anti-symplectique $\beta$ d\'efinie sur l'espace hamiltonien $M$ consid\'er\'e. Le lieu des points fixes $M^{\beta}:=Fix(\beta)$ (suppos\'e non vide) d'une telle involution est en particulier une sous-vari\'et\'e lagrangienne de $M$. En r\'ef\'erence \`a une id\'ee d'Atiyah expos\'ee dans \cite{Du} (o\`u Duistermaat \'etudie le cas des actions hamiltoniennes de tores en pr\'esence d'une involution anti-symplectique) consistant \`a voir les grassmaniennes r\'eelles comme sous-vari\'et\'es lagrangiennes des grassmaniennes complexes, les r\'esultats de convexit\'e concernant l'image par l'application moment de la sous-vari\'et\'e $M^{\beta}$ portent le nom de th\'eor\`emes de convexit\'e \emph{r\'eels}. Si un tel r\'esultat est connu dans le cadre hamiltonien usuel (voir \cite{OSS}), ce n'est \`a notre connaissance pas le cas dans le cadre quasi-hamiltonien, o\`u l'on dispose uniquement d'un r\'esultat de convexit\'e concernant l'image par le moment de la vari\'et\'e $M$ tout enti\`ere, d\^u \`a Meinrenken et Woodward (voir \cite{MW2} et \cite{AMM}).\\
\indent Le but de la pr\'esente note est de donner l'\'enonc\'e et le principe de la preuve d'un tel th\'eor\`eme de convexit\'e r\'eel pour les actions quasi-hamiltoniennes de groupes de Lie compact connexes et simplement connexes (th\'eor\`eme \ref{convexite_max}). 

\section{Enonc\'e du th\'eor\`eme}

Le th\'eor\`eme suivant est d\'emontr\'e en d\'etail dans \cite{these_Sch}. Nous indiquons dans la section \ref{preuve} les principales id\'ees et \'etapes de la d\'emonstration. Dans toute la suite, on se donne un groupe de Lie compact connexe et simplement connexe $G$, dont l'alg\`ebre de Lie $\mathfrak{g}=Lie(G)$ est suppos\'ee munie d'un produit scalaire (d\'efini positif) $Ad$-invariant $(.\ |\ .)$, et un $G$-espace quasi-hamiltonien connexe $(M,\w,\mu:M\to G)$ (notion introduite dans \cite{AMM}). Pour pouvoir \'etudier les propri\'et\'es de convexit\'e de l'application moment $\mu$, on se donne un domaine fondamental de l'action de $G$ sur lui-m\^eme par conjugaison. Un tel domaine est de la forme $\exp(\clWalc)$ o\`u $\Walc$ est une \emph{alc\^ove de Weyl} de l'alg\`ebre de Lie $\mathfrak{g}=Lie(G)$ (voir \cite{Bourbaki} ou \cite{Lo}). Pour obtenir un tel domaine fondamental, on fixe un tore maximal $T\subset G$ du groupe compact connexe $G$, et on note $R$ le syst\`eme de racines associ\'e; le choix d'une chambre de Weyl $\mathfrak{t}_+ \subset \mathfrak{t}:=Lie(T)$ (rappelons que le produit scalaire $Ad$-invariant de $\mathfrak{g}$ permet d'identifier les alg\`ebres de Lie consid\'er\'es \`a leurs duaux respectifs), ou de mani\`ere \'equivalente celui d'un systeme de racines positives $R^+\subset R$, d\'etermine alors une unique alc\^ove de Weyl $\Walc\subset \mathfrak{t}\subset \mathfrak{g}$ dont l'adh\'erence contient $0\in\mathfrak{t}$. L'application exponentielle $\exp:\mathfrak{g}\to G$ induit un hom\'eomorphisme du polytope convexe compact $\clWalc$ sur le sous-ensemble ferm\'e $\exp(\clWalc)\subset G$ et une partie $A\subset \exp(\clWalc)$ est dite convexe si elle s'identifie, via l'hom\'eomorphisme $\exp:\clWalc \to \exp(\clWalc)$, \`a une partie convexe de $\clWalc\subset \mathfrak{g}$.\\
\indent Pour obtenir un th\'eor\`eme de convexit\'e \emph{r\'eel}, on suppose le groupe $G$ muni d'un automorphisme involutif $\tau$ et on note $\taum$ l'involution de $G$ d\'efinie par $\taum(g):=\tau(g^{-1})$. L'hypoth\`ese cruciale pour le th\'eor\`eme que nous pr\'esentons dans cette note consiste \`a supposer que l'involution $\taum$ v\'erifie la propri\'et\'e (P) suivante~:
\begin{itemize}
\item[\textbf{(P)}] Il existe un tore \emph{maximal} $T$ de $G$ fix\'e \emph{point par point} par l'involution $\taum$.\label{hyp_P}
\end{itemize}
Un tore de $G$ fix\'e point par point par $\taum$ existe toujours mais ce n'est pas n\'ecessairement un tore \emph{maximal} de $G$. Dans ce dernier cas, la paire sym\'etrique $(G,\tau)$ est dite \emph{de rang maximal} dans \cite{Lo}. Une involution $\tau$ telle que $\taum$ v\'erifie la propri\'et\'e (P) existe toujours sur un groupe de Lie compact connexe simplement connexe $G$ donn\'e (m\^eme sans l'hypoth\`ese de simple connexit\'e, voir \cite{Lo}, et pour un exemple consid\'erer $SU(n)$ ou $U(n)$ muni de l'involution $\tau(u)=\overline{u}$). Si $T$ est un tore maximal de $G$ fix\'e point par point par point par $\taum$, on a en particulier $\exp(\clWalc)\subset T\subset Fix(\taum)$. On a alors le th\'eor\`eme suivant~:

\begin{theoreme}[Un th\'eor\`eme de convexit\'e r\'eel dans le cadre quasi-hamiltonien]\label{convexite_max}
Soit $(G,\tau)$ un groupe de Lie compact connexe et simplement connexe muni
d'un automorphisme involutif $\tau$ tel que l'involution 
$\taum:g\mapsto \tau(g^{-1})$ laisse fixe point par point un tore maximal $T$ de $G$
et soit $\Walc\subset \mathfrak{t}:=Lie(T)$ une alc\^ove de Weyl.
Soit $(M,\w,\mu:M\to G)$ un $G$-espace quasi-hamiltonien connexe dont
l'application moment $\mu:M\to G$ est propre et soit $\beta:M\to M$ une
involution de $M$ v\'erifiant~:
\begin{enumerate}
\item[(i)] $\beta^*\w=-\w$
\item[(ii)] $\beta(g.x)=\tau(g).\beta(x)$ pour tout $x\in M$ et tout $g\in G$
\item[(iii)] $\mu\circ\beta=\taum\circ\mu$
\item[(iv)] $M^{\beta}:=Fix(\beta)\not=\emptyset$
\item[(v)] $\mu(M^{\beta})$ a une intersection non vide avec la composante connexe $Q_0$ de $1$ dans $Fix(\taum)\subset G$
\end{enumerate}
\noindent Alors: $$\mu(M^{\beta})\cap\exp(\clWalc) =\mu(M)\cap\exp(\clWalc)$$
En particulier, $\mu(M^{\beta})\cap\exp(\clWalc)$ est un sous-polytope
convexe de $\exp(\clWalc)\simeq\clWalc\subset\mathfrak{t}$, \'egal
au polytope moment $\mu(M)\cap\exp(\clWalc)$ tout entier.
\end{theoreme}

L'hypoth\`ese (v), de nature topologique, nous est impos\'ee par le fait que, m\^eme lorsque $G$ est simplement connexe, $Fix(\taum)$ n'est pas n\'ecessairement connexe (voir \cite{Lo}). Ce probl\`eme ne se pr\'esente pas dans le cas hamiltonien usuel trait\'e dans \cite{OSS} car $Fix(\taum)\subset\mathfrak{g}^*$ est alors un espace vectoriel et est donc connexe.

\section{Principe de la preuve}\label{preuve}

\subsection{Existence d'une tranche symplectique}

La preuve du th\'eor\`eme \ref{convexite_max} repose sur l'existence d'une \emph{tranche symplectique} $N\subset M$ qui permet de ramener l'\'etude de $\mu(M)\cap\exp(\clWalc)$ \`a celle de $\mu(N)$, o\`u $N$ est un espace \emph{hamiltonien usuel} pour l'action du \emph{tore} maximal $T\subset G$ fix\'e dans l'\'enonc\'e du th\'eor\`eme \ref{convexite_max}. Cette id\'ee, d\^ue \`a Guillemin et Sternberg dans \cite{GS1}, est reprise dans \cite{CDM} et \cite{HNP}~: elle permet, dans le cas hamiltonien usuel, de prouver le th\'eor\`eme de Kirwan en se ramenant au th\'eor\`eme de Atiyah-Guillemin-Sternberg. C'est cette m\^eme id\'ee qui permet \`a Meinrenken et Woodward de montrer la convexit\'e de $\mu(M)\cap\exp(\clWalc)$ dans le cas quasi-hamiltonien (voir \cite{AMM}). Nous poursuivons ici cette approche g\'eom\'etrique pour obtenir la convexit\'e de $\mu(M^{\beta})\cap\exp(\clWalc)$. Notons que la d\'emonstration de O'Shea et Sjamaar d'un th\'eor\`eme de convexit\'e r\'eel dans le cadre hamiltonien usuel repose sur id\'ees diff\'erentes (voir notamment \cite{Sjamaar}), si bien que la strat\'egie que nous suivons ici sugg\`ere la possibilit\'e d'obtenir une preuve diff\'erente du th\'eor\`eme de O'Shea et Sjamaar. Pr\'ecisons ce que signifie l'existence d'une tranche symplectique~:
\begin{proposition}[Existence d'une tranche symplectique, \cite{AMM}]\label{tranche}
Soit $G$ un groupe de Lie compact connexe et simplement connexe et soit $(M,\w,\mu:M\to G)$ un $G$-espace quasi-hamiltonien. On suppose que l'application moment $\mu:M\to G$ est propre. Soit $T\subset G$ un tore maximal de $G$ et soit $\Walc\subset \mathfrak{t}:=Lie(T)$ une alc\^ove de Weyl. On note $p: G\to G/Int(G)$ la projection de $G$ sur l'espace des classes de conjugaison et on rappelle que l'exponentielle $\exp:\mathfrak{t}\to T$ induit un hom\'eomorphisme $\clWalc\simeq G/Int(G)$ (c'est-\`a-dire que $\exp(\clWalc)\subset G$ est un domaine fondamental pour l'action par conjugaison).\\ Alors il existe une sous-vari\'et\'e $N\subset M$ telle que~:
\begin{enumerate}
\item[(i)] $N$ est connexe
\item[(ii)] $N$ est stable par l'action de $T$
\item[(iii)] la $2$-forme $\w|_N$ est symplectique
\item[(iii)] l'action de $T$ sur $(N,\w|_N)$ est hamiltonienne avec pour application moment $$\tilde{\mu}:= p\circ\mu|_N : N\longrightarrow G/Int(G)\simeq \clWalc \subset \mathfrak{t}$$
\item[(v)] l'ensemble $G.N:=\{g.x~:~x\in N, g\in G\}$ est dense dans $M$ et l'ensemble $\mu(N)$ est dense dans $\mu(M)\cap\exp(\clWalc)$
\end{enumerate}
\end{proposition}

La sous-vari\'et\'e $N$ dont l'existence est garantie par la proposition \ref{tranche} est appel\'ee une \emph{tranche symplectique} car c'est une vari\'et\'e symplectique v\'erifiant $\overline{G.N}=M$. Ainsi qu'indiqu\'e, ce r\'esultat a sa source dans \cite{AMM}. Cependant, la preuve qui en est est donn\'ee dans \cite{these_Sch} est de nature diff\'erente~: on y suit la d\'emarche de \cite{HNP} pour montrer qu'il est possible d'obtenir une telle tranche symplectique $N\subset M$ qui soit de la forme $N=\mu^{-1}(\exp(\Walc_S))$ o\`u $\Walc_S$ est une \emph{face} du polytope convexe $\clWalc$. Pr\'ecis\'ement, la face $\Walc_S$ est celle qui contient les points de $\mu(M)$ dont la classe de conjugaison est de dimension la plus grande possible (on montre en particulier qu'une telle face est unique).\\
\indent Le fait que $\mu(N)$ soit dense dans $\mu(M)\cap\exp(\clWalc)$ montre qu'il suffit, pour montrer la convexit\'e de $\mu(M)\cap\exp(\clWalc)$, de montrer celle de $\mu(N)$. Puisque, d'apr\`es la proposition \ref{tranche}, $N$ est un espace hamiltonien au sens usuel pour l'action du \emph{tore} $T$, on est ramen\'e au th\'eor\`eme d'Atiyah-Guillemin-Sternberg. En r\'ealit\'e, la preuve est techniquement plus subtile car l'application moment $\mu|_N$ n'est pas n\'ecessairement propre, mais on peut lever ces difficult\'es de nature topologique gr\^ace \`a l'application du \emph{principe local-global} (voir \cite{HNP}), que nous ne d\'etaillerons pas ici (voir aussi \cite{Benoist}).

\subsection{Image du lieu des points fixes d'une involution}

Pour montrer la convexit\'e de $\mu(M^{\beta})\cap\exp(\clWalc)$ (et d\'ecrire cet ensemble comme sous-polytope du polytope moment $\mu(M)\cap\exp(\clWalc)$), l'id\'ee, dans la logique de la construction ci-dessus, est d'appliquer le th\'eor\`eme de Duistermaat \`a la tranche symplectique $N\subset M$ obtenue dans la proposition \ref{tranche}. La premi\`ere difficult\'e est de montrer que $N^{\beta}\not=\emptyset$. C'est ici qu'intervient l'hypoth\`ese que $\exp(\clWalc)\subset T \subset Fix(\taum)$. En effet, $N$ \'etant de la forme $\mu^{-1}(\exp(\Walc_S))$ avec $\Walc_S\subset \clWalc$, la relation $\mu\circ\beta=\taum\circ\mu$ du th\'eor\`eme \ref{convexite_max}, jointe au faite que $T$ est fix\'e \emph{point par point} par $\taum$, implique $\mu\circ\beta=\mu$, ce qui est l'une des hypoth\`eses du th\'eor\`eme de Duistermaat (voir \cite{Du}) et montre en particulier que $\beta(N)\subset N$. Il reste alors \`a montrer que $N^{\beta}$ est non-vide, ce qui n'est pas imm\'ediat. Plus pr\'ecis\'ement, l'\'etape cruciale de notre raisonnement est le r\'esultat suivant:

\begin{lemme}\label{cle}
Soit $(M,\w,\mu:M\to G)$ un $G$-espace quasi-hamiltonien muni d'une involution $\beta$ v\'erifiant les hypoth\`eses du th\'eor\`eme \ref{convexite_max} et soit $N\subset M$ une tranche symplectique comme dans la proposition \ref{tranche}. D'apr\`es l'hypoth\`ese (v) du th\'eor\`eme \ref{convexite_max}, il existe $x_0\in M^{\beta}$ tel que $\mu(x_0)$ soit contenu dans la composante connexe $Q_0$ de $1$ dans $Fix(\taum)\subset G$. Notons $L_0$ la composante connexe de $x_0$ dans $M^{\beta}$. Remarquons que l'on a \'egalement $\exp(\clWalc)\subset Q_0\subset Fix(\taum)$. Alors, on a~:
\begin{enumerate}
\item[(i)] $N$ est stable par $\beta$ et $N^{\beta}\not=\emptyset$
\item[(ii)] $\mu(N^{\beta})$ est dense dans $\mu(L_0)\cap\exp(\Walc)$
\end{enumerate}
\end{lemme}

Comme pour montrer la convexit\'e de $\mu(N)$, il faut ensuite appliquer le th\'eor\`eme de Duistermaat de mani\`ere appropri\'ee pour montrer que $\mu(N^{\beta})=\mu(N)$. On d\'eduit alors de la proposition \ref{tranche} et du lemme \ref{cle} que~:
$$\mu(M^{\beta})\cap\exp(\clWalc)\supset\mu(L_0)\cap\exp(\clWalc) = \overline{\mu(N^{\beta})}=\overline{\mu(N)} = \mu(M)\cap\exp(\clWalc)$$ ce qui prouve le th\'eor\`eme \ref{convexite_max}, l'inclusion r\'eciproque \'etant \'evidente. La preuve du fait que $N^{\beta}\not=\emptyset$ repose  sur le fait que $N$ est de la forme $\mu^{-1}(\exp(\Walc_S))$ o\`u $\Walc_S$ est une face du polytope convexe $\clWalc\subset \mathfrak{t}$ et sur le fait que $\exp(\clWalc)\subset Fix(\taum)$. Sans cette derni\`ere hypoth\`ese, qui est notre hypoth\`ese (P) page \pageref{hyp_P}, rien ne garantit que $exp(\Walc_S)$ contienne des points de $\mu(M^{\beta})$~: la relation $\mu\circ\beta=\taum\circ\mu$ du th\'eor\`eme \ref{convexite_max} entra\^ine que $\mu(M^{\beta})\subset Fix(\taum)$ mais si le domaine fondamental $\exp(\clWalc)$ n'est pas fix\'e tout entier par $\taum$ on peut par exemple avoir, selon la dimension de la face $\exp(\Walc_S)$, une intersection vide $\exp(\Walc_S)\cap Fix(\taum)=\emptyset$ avec le lieu des points fixes de $\taum$.\\
\indent Un exemple d'involution $\beta$ v\'erifiant les hypoth\`eses du th\'eor\`eme \ref{convexite_max} est donn\'e dans \cite{these_Sch}. En guise d'application, on montre comment d\'eduire de ce r\'esultat l'existence d'une sous-vari\'et\'e lagrangienne remarquable dans l'espace des repr\'esentations $\mathrm{Hom}_{\mathcal{C}}(\pi,G)/G$ du groupe fondamental $\pi:=\pi_1(S^2\backslash\{ s_1,...,s_l\})$ dans le groupe de Lie compact connexe et simplement connexe $G$ fix\'e initialement.

\end{document}